\documentstyle[11pt]{article}
\title{On the weak Deligne-Simpson problem for index of rigidity 2
\footnote{Research partially supported 
by INTAS grant 97-1644}}
\author{Vladimir Petrov Kostov\\ \\ \hspace{7cm}
{\sl To the memory of my mother}\\ \\  Appendix by Ofer Gabber} 
\date{}
\newtheorem{tm}{Theorem}
\newtheorem{lm}[tm]{Lemma}
\newtheorem{cor}[tm]{Corollary}
\newtheorem{prop}[tm]{Proposition}
\newtheorem{rem}[tm]{Remark}
\newtheorem{rems}[tm]{Remarks}
\newtheorem{defi}[tm]{Definition}
\newtheorem{ex}[tm]{Example}

\newtheorem{nota}{Notation}
\begin{document}
\maketitle 

\section{Introduction}
\subsection{Regular and Fuchsian linear systems}

Consider the linear system of ordinary differential equations defined on 
Riemann's sphere:

\begin{equation}\label{system}
{\rm d}X/{\rm d}t=A(t)X
\end{equation}
where the $n\times n$-matrix $A$ is meromorphic on ${\bf C}P^1$, with poles 
at $a_1$, $\ldots$, $a_{p+1}$; the dependent variables $X$ 
form an $n\times n$-matrix. Without loss of generality we assume that 
$\infty$ is not among the poles $a_j$ and not a pole of the $1$-form 
$A(t)$d$t$.

\begin{defi}
System (\ref{system}) is called {\em regular} at the pole $a_j$ if its 
solutions have a {\em moderate growth rate} there, i.e. for every 
sector $S$ centered at $a_j$ and of sufficiently small radius and for 
every solution $X$ restricted to the sector 
there exists $N_j\in {\bf R}$ such that $||X(t-a_j)||=O(|t-a_j|^{N_j})$ 
for all $t\in S$. System (\ref{system}) is {\em regular} if it is 
regular at all poles $a_j$. 

System (\ref{system}) is {\em Fuchsian} if its poles are logarithmic. 
Every Fuchsian system is regular.
\end{defi}

A Fuchsian system admits the presentation 

\begin{equation}\label{Fuchs} 
{\rm d}X/{\rm d}t=(\sum _{j=1}^{p+1}A_j/(t-a_j))X~,~
A_j\in gl(n,{\bf C})
\end{equation}
The sum of its {\em matrices-residua} $A_j$ equals 0, i.e. 

\begin{equation}\label{A_j}
A_1+\ldots +A_{p+1}=0
\end{equation}

When one performs a linear change of the dependent variables 

\begin{equation}\label{W}
X\mapsto W(t)X
\end{equation}
$W$ being meromorphic on 
${\bf C}P^1$, then system (\ref{system}) changes as follows:

\begin{equation}\label{gauge}
A\rightarrow -W^{-1}({\rm d}W/{\rm d}t)+W^{-1}AW
\end{equation}
(i.e. the system undergoes the {\em gauge transformation}). This 
transformation preserves regularity but, in general, it does not 
preserve being Fuchsian. The only invariant under the group 
of linear transformations (\ref{gauge}) is 
the {\em monodromy group} of the system.

Set $\Sigma ={\bf C}P^1\backslash \{ a_1,\ldots ,a_{p+1}\}$. 
Fix a base point $a_0\in \Sigma$ and a matrix $B\in GL(n,{\bf C})$. 

\begin{defi}
A {\em monodromy operator} of system (\ref{system}) defined by the class 
of homotopy equivalence in $\Sigma$ of a closed contour $\gamma$ with base 
point $a_0$ and bypassing the poles of the system 
is a linear operator $M$ acting on the  solution space of the system which 
maps the solution $X$ with $X|_{t=a_0}=B$ into the value of its analytic 
continuation along $\gamma$. Notation: $X\stackrel{\gamma}{\mapsto}XM$. 

The {\em monodromy group} is the subgroup of $GL(n,{\bf C})$ generated 
by all monodromy operators. 
\end{defi}

\begin{rem}\label{anti}
The monodromy group is an {\em antirepresentation} 
$\pi _1(\Sigma )\rightarrow GL(n,{\bf C})$ because one has 

\[ X\stackrel{\gamma _1}{\mapsto}XM_1\stackrel{\gamma _2}{\mapsto}
XM_2M_1~~~~~~(*),\] 
i.e. the concatenation $\gamma _1\gamma _2$ of the two contours defines the 
monodromy operator $M_2M_1$.   
\end{rem}

One usually chooses a standard set of contours $\gamma _j$, $j=1, \ldots ,p+1$ 
defining the generators $M_j$ of the monodromy group as follows. 
One connects $a_0$ with the points $a_j'$ (where $a_j'$ is close to $a_j$) by 
simple Jordan curves $\delta _j$ which intersect two by two only at $a_0$. The 
contour $\gamma _j$ consists of $\delta _j$, of a small circumference centered 
at $a_j$ and passing through $a_j'$ (run counterclockwise) and of 
$\delta _j$ run 
from $a_j'$ to $a_0$. Thus $\gamma _j$ is freely homotopic to a small loop 
circumventing counterclockwise $a_j$ (and no other pole $a_i$). 
The indices of the poles are 
chosen such that the indices of the contours increase from 1 to $p+1$ 
when one turns around $a_0$ clockwise.  

For the standard choice of the contours the generators $M_j$ satisfy the 
relation 

\begin{equation}\label{M_j}
M_1\ldots M_{p+1}=I
\end{equation}
which can be thought of as a multiplicative analog of (\ref{A_j}) if the 
system is Fuchsian. Equality (\ref{M_j}) results from $(*)$ 
(see Remark~\ref{anti}) -- the concatenation of 
contours $\gamma _{p+1}\ldots \gamma _1$ is homotopy equivalent to 0.

\begin{rems}\label{MexpA}
1) The monodromy group is correctly defined only up to conjugacy due to 
the freedom to choose $a_0$ and $B$.

2) For a Fuchsian system the generator $M_j$ defined as above is conjugate to 
$\exp (2\pi iA_j)$ if $A_j$ has no eigenvalues differing by a non-zero 
integer.

3) The generators $M_j$ of the monodromy group when defined after a 
standard set of contours $\gamma _j$, are conjugate to the corresponding 
operators $L_j$ of {\em local monodromy}, i.e. when the poles $a_j$ are 
circumvented counterclockwise along small loops. The operators $L_j$ of 
a regular system can be 
computed (up to conjugacy) algorithmically -- one first makes the system 
Fuchsian at $a_j$ by means of a change (\ref{W}) 
as explained in \cite{Mo} and then carries out the computation 
as explained in \cite{Wa}. 
\end{rems}

\subsection{The Deligne-Simpson problem and its weak version}

A natural question to ask is whether for given 
local monodromies (around the poles 
$a_j$) defined up to conjugacy there exists a Fuchsian system with such local 
monodromies; this is a realization problem. The difficulty is that one must 
have (\ref{M_j}). A similar question can be asked for matrices $A_j$ whose 
sum is 0 (see (\ref{A_j})). The question can be made more precise:

{\em Give necessary and sufficient conditions on the choice of 
the conjugacy 
classes $C_j\subset GL(n,{\bf C})$ or $c_j\subset gl(n,{\bf C})$ so that 
there exist 
irreducible $(p+1)$-tuples of matrices $M_j\in C_j$ or $A_j\in c_j$ 
satisfying respectively (\ref{M_j}) or (\ref{A_j}).}

This is the {\em Deligne-Simpson problem (DSP)}. ``Irreducible'' 
means ``with no common proper invariant subspace''. In technical terms, 
impossible to bring the $(p+1)$-tuple to a block upper-triangular form 
with the same sizes of the diagonal blocks for all matrices $M_j$ or $A_j$ 
by simultaneous conjugation.

\begin{rems}\label{multadd}
1) A priori one does not 
expect the problem 
to include the requirement of 
irreducibility. However, for almost all possible eigenvalues of the 
conjugacy classes the monodromy group is indeed irreducible 
and for such eigenvalues (called {\em generic}, see the next definition) 
the answer to the problem depends actually not 
on the conjugacy classes but only on the Jordan normal forms which they 
define, see Theorem~\ref{generic}. 

2) In the multiplicative version (i.e. for matrices $M_j$) 
the DSP was stated by P.Deligne (in the 
additive, i.e. for matrices $A_j$, presumably by the author) and C.Simpson 
was the first to obtain 
important results towards its resolution, see \cite{Si}. The multiplicative 
version is more important because the monodromy group is invariant under 
the action of the group of linear changes (\ref{W}) while the 
matrices-residua of a Fuchsian system are not, see rule (\ref{gauge}) and the 
lines following it. The additive version is technically easier to deal with 
and one can deduce corollaries about the multiplicative one due to 2) of 
Remarks~\ref{MexpA}.
\end{rems}

We consider only such conjugacy classes $C_j$ (resp. $c_j$) for which 
the necessary condition $\prod \det (C_j)=1$ (resp. 
$\sum$Tr$(c_j)=0$) holds. (These conditions result from (\ref{M_j}) and 
(\ref{A_j}) respectively.) In terms of the eigenvalues $\sigma _{k,j}$  
(resp. $\lambda _{k,j}$) of the matrices from $C_j$ (resp. $c_j$) repeated 
with their multiplicities, these conditions read     

\begin{equation}\label{evs}
\prod _{k=1}^n\prod _{j=1}^{p+1}\sigma _{k,j}=1~~,~~{\rm resp.}~~ 
\sum _{k=1}^n\sum _{j=1}^{p+1}\lambda _{k,j}=0
\end{equation} 

\begin{defi}\label{nongenrel}
An equality of the form 

\begin{equation}\label{relation}
\prod _{j=1}^{p+1}\prod _{k\in \Phi _j}\sigma _{k,j}=1~~,~~{\rm resp.}~~ 
\sum _{j=1}^{p+1}\sum _{k\in \Phi _j}\lambda _{k,j}=0
\end{equation} 
is called a {\em non-genericity relation};  
the non-empty sets $\Phi _j$ contain one and the same number $<n$ of indices  
for all $j$. Eigenvalues satisfying none of these relations are called 
{\em generic}. Reducible  
$(p+1)$-tuples exist only for non-generic eigenvalues (the 
eigenvalues of each diagonal block of a block upper-triangular  
$(p+1)$-tuple satisfy some non-genericity relation). 
\end{defi}

\begin{rem}\label{tricky}
In the case of matrices 
$A_j$, if the greatest common divisor $d$ of the 
multiplicities of all 
eigenvalues of all $p+1$ matrices is $>1$, then a non-genericity relation 
results automatically from $\sum$Tr$(c_j)=0$. In  
the case of matrices $M_j$  
the equality $\prod \sigma _{k,j}=1$ implies that 
if one divides by $d$ the multiplicities of all eigenvalues, then their 
product would equal $\exp (2\pi ik/d)$, $0\leq k\leq d-1$, 
not necessarily 1, and 
a non-genericity relation might or might not hold (see Example~\ref{ex1}).
\end{rem} 

\begin{defi}
Call {\em Jordan normal form (JNF) of size $n$} a family 
$J^n=\{ b_{i,l}\}$ ($i\in I_l$, $I_l=\{ 1,\ldots ,s_l\}$, $l\in L$) of 
positive integers $b_{i,l}$ 
whose sum is $n$. Here $L$ is the set of indices of 
eigenvalues (all distinct) and 
$I_l$ is the set of Jordan blocks with the $l$-th eigenvalue, $b_{i,l}$ is the 
size of the $i$-th block with this eigenvalue. E.g. the JNF 
$\{ \{ 2,1\} \{ 4,3,1\} \}$ is of size $11$ and with two eigenvalues to the 
first (resp. second) of 
which there correspond two (resp. three) Jordan blocks, of sizes $2$ and $1$ 
(resp. $4$, $3$ and $1$). An $n\times n$-matrix 
$Y$ has the JNF $J^n$ (notation: $J(Y)=J^n$) if to its distinct 
eigenvalues $\lambda _l$, $l\in L$, there belong Jordan blocks of sizes 
$b_{i,l}$. We denote by $J(C)$ the JNF defined by the conjugacy class $C$. 
\end{defi}

\begin{ex}\label{ex1} 
Let $p+1=n=4$, let for each $j$ the JNF $J_j^4$ consist 
of two Jordan blocks $2\times 2$, with equal eigenvalues. If the eigenvalues 
of $M_1$, $\ldots$, $M_4$ equal $i$, $1$, $1$ and $1$ 
(resp. $-1$, $1$, $1$, $1$), then 
they are generic (resp. they are non-generic) -- when their multiplicities are 
reduced twice, then their product equals $-1$, a primitive root of $1$ of 
order $2$ (resp. their product equals $1$).
\end{ex}

When the eigenvalues are not generic, then the $(p+1)$-tuples 
of matrices $M_j\in C_j$ or $A_j\in c_j$ (if it exists) is not automatically  
irreducible. Therefore we give the following

\begin{defi}
The formulation of the 
{\em weak Deligne-Simpson problem} is obtained when in the one 
of the DSP one replaces the requirement of irreducibility by 
the weaker requirement the centralizer of the $(p+1)$-tuple of matrices 
$A_j$ or $M_j$ to be trivial, i.e. reduced to scalars.
\end{defi}

\begin{rem}
If one states the problem asking only the matrices 
$M_j\in C_j$ or $A_j\in c_j$ to 
satisfy respectively condition (\ref{M_j}) or (\ref{A_j}) and with no 
requirement of irreducibility or triviality of the centralizer, then solving 
the problem becomes much harder and the answer to it depends essentially 
on the eigenvalues (not only on the JNFs). E.g., suppose that  
$p=n=2$ and that two of the matrices $M_j$ (resp. $A_j$) have distinct 
eigenvalues $\sigma _{1,j}\neq \sigma _{2,j}$, $j=1,2$ 
(resp. $\lambda _{1,j}\neq \lambda _{2,j}$) while the third must be scalar 
(i.e. $\sigma _{1,3}=\sigma _{2,3}$, 
resp. $\lambda _{1,3}=\lambda _{2,3}$). Then such triples 
exist exactly if $\sigma _{1,1}\sigma _{1,2}\sigma _{1,3}=1$ or 
$\sigma _{1,1}\sigma _{2,2}\sigma _{1,3}=1$ 
(resp. $\lambda _{1,1}+\lambda _{1,2}+\lambda _{1,3}=0$ 
or $\lambda _{1,1}+\lambda _{2,2}+\lambda _{1,3}=0$). Hence, such triples 
exist exactly if the eigenvalues are not generic. 
\end{rem}

\begin{defi}
Denote by $\{ J_j^n\}$  
a $(p+1)$-tuple of JNFs, $j=1$,$\ldots$, $p+1$. 
We say that the DSP (resp. the weak DSP) is {\em solvable}  for a 
given $\{ J_j^n\}$ and given eigenvalues if there exists an 
irreducible $(p+1)$-tuple (resp. a $(p+1)$-tuple with trivial centralizer) 
of matrices $M_j$ satisfying (\ref{M_j}) or of matrices $A_j$ satisfying 
(\ref{A_j}), with $J(M_j)=J_j^n$ or $J(A_j)=J_j^n$ and with the given 
eigenvalues. By definition, the DSP is solvable for 
$n=1$. Solvability of the DSP implies the one of the weak DSP.
\end{defi}

For generic eigenvalues the DSP is solved -- the result is formulated 
in \cite{Ko3} and proved in \cite{Ko4} and \cite{Ko2}. 
In the next subsection we recall this result (Theorem~\ref{generic}). 
The result is a necessary and sufficient condition on the JNFs $J(C_j)$ 
or $J(c_j)$. 

The aim of the present paper is to show an example of a large class 
of $(p+1)$-tuples of conjugacy classes $C_j$ or $c_j$ 
(such that the conditions of 
Theorem~\ref{generic} hold for the JNFs $J(C_j)$ or $J(c_j)$), with 
non-generic eigenvalues, for which the weak DSP is not solvable.

\subsection{The results known up to now}

\begin{nota}
For a conjugacy class $C$ in $GL(n,{\bf C})$ or $gl(n,{\bf C})$ denote by 
$d(C)$ its dimension and for a matrix $Y$ from $C$ set 
$r(C):=\min _{\lambda \in {\bf C}}{\rm rank}(Y-\lambda I)$. The integer 
$n-r(C)$ is the maximal number of Jordan blocks of $J(Y)$ with one and the 
same eigenvalue. Set $d_j:=d(C_j)$ (resp. $d(c_j)$), $r_j:=r(C_j)$ 
(resp. $r(c_j)$). The quantities 
$r(C)$ and $d(C)$ depend only on the JNF $J(Y)=J^n$, not 
on the eigenvalues, so we write sometimes $r(J^n)$ and $d(J^n)$.
\end{nota} 

\begin{prop}\label{d_jr_j}
(C. Simpson, see \cite{Si}.) The 
following two inequalities are necessary conditions for the solvability 
of the DSP:

\[ \begin{array}{lll}d_1+\ldots +d_{p+1}\geq 2n^2-2&~~~~~&(\alpha _n)\\
{\rm for~all~}j,~r_1+\ldots +\hat{r}_j+\ldots +r_{p+1}\geq n&~~~~~&
(\beta _n)\end{array}\]
\end{prop}

The following condition is not necessary and in most cases it is sufficient 
for the solvability of the DSP, see~\cite{Ko3} and~\cite{Ko1}:

\[ (r_1+\ldots +r_{p+1})\geq 2n~~~~~~~~~~~~~~~~(\omega _n)~~~.\]

For a given $\{ J_j^n\}$ with $n>1$, which satisfies conditions 
$(\alpha _n)$ and $(\beta _n)$ and doesn't satisfy condition 
$(\omega _n)$ set $n_1=r_1+\ldots +r_{p+1}-n$. Hence, $n_1<n$ and 
$n-n_1\leq n-r_j$. Define 
the $(p+1)$-tuple $\{ J_j^{n_1}\}$ as follows: to obtain the JNF 
$J_j^{n_1}$ 
from $J_j^n$ one chooses one of the eigenvalues of $J_j^n$ with 
greatest number $n-r_j$ of Jordan blocks, then decreases  
by 1 the sizes of the $n-n_1$ {\em smallest} Jordan blocks with this 
eigenvalue and deletes the Jordan blocks of size 0. For the above 
construction we use the notation $\Psi :\{ J_j^n\}\mapsto \{ J_j^{n_1}\}$. 

\begin{defi}
A $(p+1)$-tuple $\{ J_j^n\}$ with $n>1$ is {\em good} if 

1) it satisfies conditions $(\alpha _n)$ and $(\beta _n)$ and 

2) either $\{ J_j^n\}$ satisfies 
condition $(\omega _n)$ or the construction 
$\Psi$ iterated as long as it is defined 
stops at a $(p+1)$-tuple $\{ J_j^{n'}\}$ either with $n'=1$ or satisfying 
condition $(\omega _{n'})$.

By definition, a $(p+1)$-tuple of JNFs with $n=1$ is good.
\end{defi}

\begin{tm}\label{generic}
(see \cite{Ko3}, Theorem 8.) 
Let $n>1$. The DSP is solvable for the conjugacy classes $C_j$ or 
$c_j$ (with generic eigenvalues,  
defining the JNFs $J_j^n$) if and only if the $(p+1)$-tuple $\{ J_j^n\}$ 
is good.
\end{tm}

\begin{rem}
The quantity $\kappa =2n^2-\sum _{j=1}^{p+1}d_j$ 
(called {\em index of rigidity}, see \cite{Ka})  
is invariant for the 
construction $\Psi$, see \cite{Ko3}. Therefore one can drop 
condition $(\alpha _n)$ in the definition of a good $(p+1)$-tuple -- 
condition $(\alpha _{n'})$ always holds for the $(p+1)$-tuple of JNFs 
$\{ J_j^{n'}\}$, see \cite{Ko3} and \cite{Ko1}; if $n'=1$, 
then $(\alpha _{n'})$ is an 
equality, if there holds $(\omega _{n'})$, then $(\alpha _{n'})$ holds and 
is a strict inequality. 
\end{rem}

\subsection{The basic result\protect\label{basres}}

In the present paper we consider the case when the index of rigidity 
equals 2. 

\begin{rem}\label{notboth}
In this case if there exist irreducible $(p+1)$-tuples, then 
they are unique up to conjugacy and the coexistence of irreducible and 
reducible $(p+1)$-tuples is impossible, see \cite{Ka} and \cite{Si} for the 
multiplicative version. 

In the additive version this is also true -- if there 
exist an irreducible and a reducible $(p+1)$-tuples of matrices $A_j\in c_j$, 
then one can multiply them by a non-zero scalar $h$ so that there should be no 
non-zero integer differences between any two eigenvalues of a given matrix 
$hA_j$ and any integer sum of eigenvalues of the matrices should be $0$. The 
monodromy operators of a Fuchsian system with matrices-residua $hA_j$ are 
conjugate to $\exp (2\pi ihA_j)$. The monodromy group of the system with 
(ir)reducible $(p+1)$-tuple of matrices-residua $hA_j$ is (ir)reducible as 
well; in the reducible case this is evident, in the irreducible one 
this follows from Theorem~5.1.2 from \cite{Bo} (the latter states that if the 
monodromy group is reducible and if the 
sum of the exponents relative to an invariant subspace is $0$, then the 
matrices-residua can be simultaneously conjugated to a block 
upper-triangular form; the sum of the exponents for the system with residua 
$hA_j$ is a sum of eigenvalues of these matrices; thus the irreducibility of 
the $(p+1)$-tuple of matrices-residua implies the one of the monodromy group). 
This is a contradiction with the non-coexistence 
of irreducible and reducible $(p+1)$-tuples in the multiplicative case.
\end{rem} 

\begin{defi}\label{subordinate}
We say that the conjugacy class $c'$ (in $gl(n,{\bf C})$ 
or in $GL(n,{\bf C})$) is {\em subordinate} to the conjugacy 
class $c$ if $c'$ belongs to the closure of $c$. This means that the  
eigenvalues of $c$ and $c'$ are the same and of the same multiplicities 
and for each eigenvalue $\lambda _i$ and for each $j\in {\bf N}$ one has 
rk$(A-\lambda _iI)^j\geq$rk$(A'-\lambda _iI)^j$ for $A\in c$, $A'\in c'$. 
If $c'\neq c$, then at least one inequality is strict.
\end{defi}

\begin{defi}\label{special}
Let $n=ln_1$, $l,n_1\in {\bf N}^*$, $n_1>1$. 
The $(p+1)$-tuple of conjugacy classes $C_j$ or $c_j$ with $\kappa =2$ 
is called $l$-{\em special} if 
for each class $C_j$ (or $c_j$) there exists a class $C_j'$ (or $c_j'$) 
subordinate to it which is a direct sum of $n_1$ copies of a conjugacy 
class $C_j''\subset GL(l,{\bf C})$ (or $c_j''\subset gl(l,{\bf C})$) 
where the $(p+1)$-tuple of JNFs $J(C_j'')$ (or $J(c_j'')$) is good and the 
product of the eigenvalues of the classes $C_j''$ equals 1 
(see Remark~\ref{tricky} and Example~\ref{ex1}; for the classes $c_j''$ 
the sum of their eigenvalues is automatically $0$). 
\end{defi}

\begin{rems}\label{c_jrig}
1) The index of rigidity of the $(p+1)$-tuple of conjugacy classes $c_j''$ 
or $C_j''$ equals 2. Indeed, one has $d(c_j)\geq d(c_j')=(n_1)^2d(c_j'')$ 
and if 
$\sum _{j=1}^{p+1}d(c_j'')\geq 2l^2$, then $\sum _{j=1}^{p+1}d(c_j)\geq 2n^2$, 
i.e. the index of rigidity of the $(p+1)$-tuple of conjugacy classes $c_j$ 
must be non-positive -- a contradiction. The reasoning holds in the 
case of classes $C_j$ as well. 

2) It follows from the above definition that in the case of matrices $A_j$ the 
eigenvalues of an $l$-special 
$(p+1)$-tuple of JNFs cannot be generic -- their multiplicities are divisible 
by $n_1$ and, hence, they satisfy a non-genericity relation, 
see Remark~\ref{tricky}. Notice that in the case of matrices $M_j$ 
the divisibility by $n_1$ alone of the multiplicities does not imply that the 
eigenvalues are 
not generic, see Remark~\ref{tricky} and Example~\ref{ex1}. Therefore the 
requirement the product of the eigenvalues of the classes $C_j''$ to equal $1$ 
(see the definition) is essential.
\end{rems}

\begin{defi}
A $(p+1)$-tuple of conjugacy classes in $gl(n,{\bf C})$ or $GL(n,{\bf C})$ 
is called {\em special} if it is 
$l$-special for some $l$. If in addition for this $l$ the classes $c_j''$ or 
$C_j''$ are diagonalizable, then the $(p+1)$-tuple is called 
{\em special-diagonal}.
\end{defi}

\begin{ex}
The triple of conjugacy classes $C_1$, $C_2$, 
$C_3\subset GL(4,{\bf C})$ (or   
$c_1,c_2,c_3\subset gl(4,{\bf C})$) defining the JNFs $J_1^4=\{ \{ 4\} \}$, 
$J_2^4=\{ \{ 1,1\} ,\{ 2\} \}$, $J_3^4=\{ \{ 1,1\} ,\{ 1,1\} \}$, 
is good (to be checked directly). 

This triple is also $2$-special -- one can choose as subordinate classes 
ones in which $C_1'$ (resp. $c_1'$) has two Jordan blocks of size 2 while 
$C_2'$ and $C_3'$ (resp. $c_2'$ and $c_3'$) define the JNF $J_3^4$. 
For such a choice 
the triple $(C_1,C_2,C_3)$ (resp. $(c_1,c_2,c_3)$) will be 2-special and 
the conjugacy classes $C_j''\subset GL(2,{\bf C})$ (resp. 
$c_j''\subset gl(2,{\bf C})$ have two distinct eigenvalues for $j=2,3$ and 
one Jordan block of size 2 for $j=1$.
\end{ex}

\begin{ex}
For $n>1$ a good $(p+1)$-tuple of unipotent conjugacy classes in 
$GL(n,{\bf C})$ 
or of nilpotent conjugacy classes in $gl(n,{\bf C})$ is 1-special, hence, 
it is special.
\end{ex}  

\begin{ex}
For $n=9$ the triple of conjugacy classes $c_j$ 
defining the JNFs $\{ \{ 2,2,1,1\} ,\{ 1,1,1\} \}$ for $j=1,2$ and 
$\{ \{ 2,2,1,1\} ,\{ 2,1\} \}$ for $j=3$ is good. Although the multiplicities 
of all eigenvalues are divisible by 3, the triple is not special (a priori 
if it is special, then it is 3-special). Indeed, the JNFs are such that 
the conjugacy classes $c_j''$ from the definition of a special $(p+1)$-tuple 
must be diagonalizable (for each eigenvalue 
of $c_j$ there are at most two Jordan blocks of size $>1$ and this size is 
actually 2). But then $c_j''$ must have each two eigenvalues, of 
multiplicities 1 and 2, which means that the triple $J(c_1'')$, $J(c_2'')$, 
$J(c_3'')$ is not good.
\end{ex}  

The basic result of the paper is the following 

\begin{tm}\label{basicresult}
The weak DSP is not solvable for special-diagonal $(p+1)$-tuples of 
conjugacy classes.
\end{tm}

\begin{rem}\label{goodisnecessary}
The $(p+1)$-tuple of conjugacy classes to be good is a necessary condition 
for the solvability of the weak DSP for index of rigidity 2, see the proof of 
this in Remark~\ref{goodisnecessary1}. (In fact, it is necessary for any 
index of rigidity $\leq 2$ but we do not need the proof of this statement in 
the present paper.)
\end{rem}

\begin{rem}
The theorem raises the following two natural questions:

1) whether it remains true for special (not necessarily special-diagonal) 
$(p+1)$-tuples of conjugacy classes;

2) whether for index of rigidity $2$ and for $(p+1)$-tuples of 
conjugacy classes defining good $(p+1)$-tuples of 
JNFs the weak DSP is unsolvable only when the $(p+1)$-tuples of 
conjugacy classes are special.

It seems that the answer to the first of them is positive although the author 
was unable to find a complete proof of it. The answer to the second question 
is not known to the author.
\end{rem}

The next subsection contains the plan of the proof of the theorem. The rest 
of the proof is given in Section~\ref{PBR}.

\subsection{Plan of the proof of Theorem \protect\ref{basicresult}}

We consider first the particular case when the conjugacy classes $C_j''$, 
resp. $c_j''$, are diagonalizable and with generic eigenvalues; we assume 
also that all non-genericity relations for 
the classes $c_j$ or $C_j$ are obvious ones, i.e. multiples of the fact that 
the sum of the eigenvalues of the classes $c_j''$ is $0$ or that the 
product of the eigenvalues of the classes $C_j''$ is $1$. Call this 
case and such special-diagonal $(p+1)$-tuples of conjugacy classes 
{\em quasi-generic}. The proof in the general case is deduced from the 
quasi-generic one using a method for 
deforming analytically $(p+1)$-tuples of matrices $A_j$ or $M_j$ (the 
method is called 
the {\em basic technical tool}, it is developed in Section~\ref{BTT}). Namely, 
if there exists a $(p+1)$-tuple with trivial centralizer of matrices $A_j$ 
or $M_j$ from a special-diagonal $(p+1)$-tuple of conjugacy classes, then it 
can be analytically deformed into one defining the same JNFs, for which 
the classes $c_j''$ or $C_j''$ have generic eigenvalues and which is with 
trivial centralizer. The deformation 
changes the eigenvalues but not their multiplicities. The possibility to 
deform the eigenvalues into ones from the quasi-generic case follows from the 
fact that the classes $c_j''$ or $C_j''$ are semisimple and the index of 
rigidity of their $(p+1)$-tuple is $2$; hence, the multiplicities of their 
eigenvalues are not divisible by an integer $>1$. Non-solvability of the 
weak DSP for quasi-generic special-diagonal $(p+1)$-tuples implies its 
non-solvability for any special-diagonal $(p+1)$-tuples.  


The following proposition was suggested by the author and 
proved by Ofer Gabber (see the proof in the Appendix).

\begin{prop}\label{Gabber}
Suppose that the index of rigidity of the $(p+1)$-tuple of 
conjugacy classes $C_j\subset GL(n,{\bf C})$ or $c_j\subset gl(n,{\bf C})$ 
equals 2. Suppose that the conjugacy classes $C_j^*$ (or $c_j^*$) are 
subordinate to the respective classes $C_j$ (or $c_j$), with $C_j^*\neq C_j$ 
(or $c_j^*\neq c_j$) for at least one value of $j$. Then the existence of 
matrices $M_j\in C_j^*$ satisfying condition (\ref{M_j}) (resp. of 
matrices $A_j\in c_j^*$ satisfying condition (\ref{A_j})) implies that 
the DSP is not solvable for the $(p+1)$-tuple of conjugacy 
classes $C_j$ (resp. $c_j$).
\end{prop}

As it was mentioned in Remark~\ref{notboth}, coexistence 
of irreducible and reducible $(p+1)$-tuples of matrices $A_j$ or $M_j$ 
with index of rigidity 2 is impossible. Therefore the above 
proposition implies

\begin{cor}\label{Gabberbis}
The DSP is not solvable for special $(p+1)$-tuples of conjugacy 
classes $C_j$ or $c_j$.
\end{cor}

Indeed, for such $(p+1)$-tuples one can construct block-diagonal 
$(p+1)$-tuples of matrices $M_j\in C_j'$ or $A_j\in c_j'$ satisfying 
respectively (\ref{M_j}) or (\ref{A_j}); this follows from the $(p+1)$-tuple 
of JNFs $J(C_j'')$ or $J(c_j'')$ being good.~~~~~$\Box$ 

Further we need the following

\begin{prop}\label{dimensionrig}
For index of rigidity $\kappa$ the dimension of the variety ${\cal V}$ (when 
it is not empty) 
of $(p+1)$-tuples with 
trivial centralizers of matrices $A_j\in c_j$ or $M_j\in C_j$ equals 
$n^2+1-\kappa$.
\end{prop}

The proposition is proved in Subsection~\ref{prdimensionrig}. 

\begin{prop}\label{rigidred}
A reducible $(p+1)$-tuple with 
trivial centralizer of matrices $A_j\in c_j$ or $M_j\in C_j$, where the index 
of rigidity of the  
$(p+1)$-tuple of conjugacy classes $c_j$ or $C_j$ 
equals 2, is up to conjugacy block upper-triangular with diagonal blocks 
defining irreducible representations each with index of rigidity 2.
\end{prop}

The proposition was suggested by the author and proved by O. Gabber, see 
the proof in the Appendix.
Note that the 
$(p+1)$-tuple of conjugacy classes is not supposed to be special. 
The proposition implies 

\begin{cor}\label{rrr}
All (if any) quasi-generic special-diagonal $(p+1)$-tuples of matrices 
$A_j\in c_j$ or 
$M_j\in C_j$ which solve the weak DSP are (up to conjugacy) 
block upper-triangular and their 
diagonal blocks define irreducible representations $P_i$ each 
with index of rigidity 2.
\end{cor}

\begin{lm}\label{allequal}
If a quasi-generic special-diagonal $(p+1)$-tuple of matrices 
$A_j\in c_j$ or 
$M_j\in C_j$ is (up to conjugacy) 
block upper-triangular with  
diagonal blocks defining irreducible representations each 
with index of rigidity 2, then these representations are 
equivalent and of rank $l$.
\end{lm}

The lemma is proved in Subsection~\ref{prallequal}. It implies 

\begin{cor}\label{allequalcor}
A $(p+1)$-tuple of matrices $A_j\in c_j$ or 
$M_j\in C_j$ satisfying the conditions of Lemma~\ref{allequal} is with 
non-trivial centralizer.
\end{cor}

Indeed, block-decompose the matrices 
from $gl(n,{\bf C})$ or $GL(n,{\bf C})$ into blocks $l\times l$. Hence, the 
centralizer contains the matrix having a block equal to $I$ in position 
$(1,n/l)$ and zeros elsewhere.~~~$\Box$

Hence, the weak DSP is not solvable for quasi-generic 
special-diagonal $(p+1)$-tuples 
of conjugacy classes. Indeed, if it were solvable, then it would be solved 
only by reducible 
$(p+1)$-tuples of matrices which up to conjugacy are block upper-triangular, 
with diagonal blocks defining equivalent representations of rank $l$ and of 
index of rigidity $2$, see Corollary~\ref{rrr} and Lemma~\ref{allequal}. 
By Corollary~\ref{allequalcor}, these $(p+1)$-tuples have non-trivial 
centralizers.

This proves Theorem \ref{basicresult} in the quasi-generic case.

\section{The basic technical tool\protect\label{BTT}}

\begin{defi}
Call {\em basic technical tool} the procedure described 
below whose aim is to deform analytically a given $(p+1)$-tuple of matrices 
$A_j$ or $M_j$ with trivial centralizer 
by changing their conjugacy classes in a desired way. 
\end{defi}

Set $A_j=Q_j^{-1}G_jQ_j$, 
$G_j$ being Jordan matrices. 
Look for a 
$(p+1)$-tuple of matrices $\tilde{A}_j$ (whose sum is $0$) of the form 

\[ \tilde{A}_j=(I+\varepsilon X_j(\varepsilon ))^{-1}
Q_j^{-1}
(G_j+\varepsilon V_j(\varepsilon ))Q_j
(I+\varepsilon X_j(\varepsilon ))\] 
where 
$\varepsilon \in ({\bf C},0)$ and 
$V_j(\varepsilon )$ are given 
matrices analytic 
in $\varepsilon$; they must satisfy the condition  
tr$(\sum _{j=1}^{p+1}\varepsilon V_j(\varepsilon ))
\equiv 0$; set $N_j=Q_j^{-1}V_jQ_j$. The existence of matrices 
$X_j(\varepsilon )$ is deduced from 
the triviality of the centralizer, using the following proposition 
(see its proof and the details in \cite{Ko4}):

\begin{prop}\label{[A_j,X_j]}
The centralizer of the $p$-tuple of matrices $A_j$ ($j=1,\ldots ,p$) is 
trivial if and only if the mapping 
$(sl(n,{\bf C}))^p\rightarrow sl(n,{\bf C})$, 
$(X_1,\ldots ,X_p)\mapsto \sum _{j=1}^p[A_j,X_j]$ is surjective.
\end{prop}

Notice that one has 
$\tilde{A}_j=A_j+\varepsilon [A_j,X_j(0)]+\varepsilon N_j+o(\varepsilon )$. 
The proposition assures the existence in first approximation w.r.t. 
$\varepsilon$ of the matrices $X_j$, the existence of true matrices $X_j$ 
analytic in $\varepsilon$ follows from the implicit function theorem. 

If for $\varepsilon \neq 0$ 
small enough the eigenvalues of the matrices $\tilde{A}_j$ are generic, then 
their $(p+1)$-tuple is irreducible. In a similar way one can deform 
analytically $(p+1)$-tuples depending on a multi-dimensional parameter.

Given a $(p+1)$-tuple of matrices $M_j^1$ with  
trivial centralizer and satisfying condition (\ref{M_j}), look for matrices 
$M_j$ (whose product is $I$) of the form 

\[ M_j=(I+\varepsilon X_j(\varepsilon ))^{-1}(M_j^1+
\varepsilon N_j(\varepsilon ))(I+
\varepsilon X_j(\varepsilon ))\] 
where the given matrices $N_j$ depend analytically on 
$\varepsilon \in ({\bf C},0)$ and the product of the determinants of the 
matrices $M_j$ is $1$; one looks for $X_j$ analytic in 
$\varepsilon$. (Like in the additive version one can set 
$M_j^1=Q_j^{-1}G_jQ_j$, $N_j=Q_j^{-1}V_jQ_j$.) The existence of 
such matrices $X_j$ follows again from the triviality of the centralizer, 
see \cite{Ko4}.

\begin{rem}\label{goodisnecessary1}
As it was mentioned in Remark~\ref{goodisnecessary}, for index of rigidity 
$2$ the $(p+1)$-tuple of 
conjugacy classes to be good is a necessary condition for the solvability of 
the weak DSP. Indeed, using the basic technical tool 
one can deform a $(p+1)$-tuple of matrices $A_j$ or 
$M_j$ with trivial centralizer into a nearby one (say, of matrices $A_j'$ 
or $M_j'$) defining the so-called 
{\em corresponding} diagonal JNFs (see the definition in \cite{Ko3} or in 
\cite{Ko4}); both $(p+1)$-tuples are simultaneously good or not and 
the latter $(p+1)$-tuple is also of index of rigidity $2$. 
Hence, the multiplicities of the eigenvalues of the matrices $A_j'$ or $M_j'$ 
are not divisible by an integer $>1$ and one can choose the eigenvalues to 
be generic. But then one applies Theorem~\ref{generic} and finds out that the 
$(p+1)$-tuple of initial conjugacy classes must be good.
\end{rem}

\section{Proof of Theorem \protect\ref{basicresult}\protect\label{PBR}}

\subsection{Proof of Proposition~\protect\ref{dimensionrig}
\protect\label{prdimensionrig}}

$1^0$. Consider first the case of matrices $A_j$. Without restriction one can 
assume that $c_j\subset sl(n,{\bf C})$. 
To find the dimension of the variety ${\cal V}$ one has first to 
consider the cartesian product 
$(c_1\times \ldots \times c_p)\subset (sl(n,{\bf C}))^p$.  
Define the mapping 
$\tau :(c_1\times \ldots \times c_p)\rightarrow sl(n,{\bf C})$ by 
the rule $\tau :(A_1,\ldots ,A_p)\mapsto A_{p+1}=-A_1-\ldots -A_p$ 
(recall that there holds (\ref{A_j})). 

$2^0$. The algebraic variety ${\cal V}$  
is the intersection of the two 
varieties in $c_1\times \ldots \times c_p\times sl(n,{\bf C})$: 
the graph of the mapping $\tau$ 
and $c_1\times \ldots \times c_p\times c_{p+1}$. This intersection is 
transversal which implies the smoothness of the 
variety ${\cal V}$. Transversality follows from Proposition~\ref{[A_j,X_j]} -- 
the tangent space to the conjugacy class $c_j$ at $A_j$ 
equals $\{ [A_j,X]|X\in gl(n,{\bf C})\}$. 

$3^0$. Recall that $d_j$ denotes dim$c_j$. Hence,  

\[ {\rm dim}\, {\cal V}=
(\sum _{j=1}^pd_j)-{\rm codim}_{sl(n,{\bf C})}c_{p+1}=
(\sum _{j=1}^pd_j)-[(n^2-1)-
d_{p+1}]~.\]
Hence, 
dim${\cal V}=\sum _{j=1}^{p+1}d_j-n^2+1=2n^2-\kappa -n^2+1=n^2+1-\kappa$. 

$4^0$. The only difference in the proof in the case of matrices $M_j$ 
is that the mapping    
$(A_1,\ldots ,A_p)$ $\mapsto A_{p+1}=-A_1-\ldots -A_p$
from ${\bf 2^0}$ has to be replaced by the mapping 
\[(M_1,\ldots ,M_p)\mapsto M_{p+1}=(M_1\ldots M_p)^{-1}.\]

The reader will be able to restitute the missing technical details after 
examining the more detailed description of the basic technical tool given in 
\cite{Ko4}. The proposition is proved.~~~~~$\Box$

\subsection{Proof of Lemma~\protect\ref{allequal}\protect\label{prallequal}}

$1^0$. Recall that in the quasi-generic case the conjugacy classes $c_j''$ or 
$C_j''$ are diagonalizable, the DSP is solvable for them, and the 
irreducible representation $Q$ they define is with index of rigidity 2; 
recall that $Q$ is unique up to conjugacy, see 
Remark~\ref{notboth}. We denote by $P_i$ also the 
diagonal block defining the representation $P_i$. 

A priori in the quasi-generic case 
every representation $P_i$ is of rank $lq_i$, 
$q_i\in {\bf N}^*$, and the multiplicity of every eigenvalue of the 
matrix $A_j$ (or $M_j$) restricted to the block $P_i$ equals 
$q_i$ times its multiplicity as eigenvalue of $c_j''$ (or of $C_j''$). 
This follows from the fact that the eigenvalues of the 
conjugacy classes $c_j''$ or $C_j''$ are generic. 

$2^0$. Denote by $c_j^*$ 
(resp. by $C_j^*$) the conjugacy class of the restriction of the matrix 
$A_j$ (resp. $M_j$) to the block $P_i$. 
For every $i$ the index of rigidity of the $(p+1)$-tuple of 
conjugacy classes $c_j^*$ or $C_j^*$ equals 2 (Corollary~\ref{rrr}). 
Suppose that a given block $P_i$ is of size $q_il$ with  
$q_i>1$. Hence, the conjugacy class which is $q_i$ times the conjugacy class 
$c_j''$ (resp. $q_i$ times $C_j''$) is subordinate to $c_j^*$ (resp. to 
$C_j^*$) for $j=1,\ldots ,p+1$. Indeed, the classes $q_ic_j''$ 
(resp. $q_iC_j''$) 
and $c_j^*$ (resp. $C_j^*$) have the same eigenvalues, of the same 
multiplicities and $c_j''$ (resp. $C_j''$) is diagonalizable. 

This means that the $(p+1)$-tuple of 
conjugacy classes $c_j^*$ or $C_j^*$ is quasi-generic 
$l$-special and of size $q_il$. 

$3^0$. The DSP is not solvable for 
the $(p+1)$-tuple of 
conjugacy classes $c_j^*$ or $C_j^*$ if $q_i>1$, see 
Corollary~\ref{Gabberbis}. Hence, 
one must have $q_i=1$, i.e. all diagonal blocks are of equal size. Hence, 
the conjugacy classes of the restrictions of the matrices $M_j$ (resp. $A_j$) 
to them equal $C_j''$ (resp. $c_j''$). 
Indeed, the eigenvalues and the multiplicities are the ones of the classes 
$C_j''$ (resp. $c_j''$) -- recall that we are in the quasi-generic case. The 
presence of Jordan blocks of size $>1$ would mean that the sum of the 
dimensions $d(C(M_j|_{P_i}))$ (resp. $d(C(A_j|_{P_i}))$) is $>2n^2-2$, 
i.e. the index of rigidity of the $(p+1)$-tuple of matrices $M_j|_{P_i}$ 
(resp. $A_j|_{P_i}$) is non-positive -- a contradiction with 
Proposition~\ref{rigidred}.

$4^0$. The index of rigidity of the 
$(p+1)$-tuple of conjugacy classes $C_j''$ (resp. $c_j''$) being 2, the 
diagonal blocks define equivalent representations (see 
Remark~\ref{notboth}).~~~~~$\Box$\\

{\Large Appendix (Ofer Gabber)}\\ 

{\bf Proof of Proposition~\ref{Gabber}.}\\ 

We prove the proposition in the multiplicative version. The proof for the 
additive one can be deduced by means of a reasoning completely analogous to  
the one from Remark~\ref{notboth} and we leave it for the reader. 

Consider distinct points $a_0$, $a_1$, $\ldots$, $a_{p+1}$ in 
${\bf P}_{\bf C}^1$. Set 
$U={\bf P}_{\bf C}^1\backslash \{ a_1,\ldots,a_{p+1}\}$ and fix usual 
generators $\gamma _i$ of $\pi _1(U,a_0)$. Set 
$j:U\hookrightarrow {\bf P}_{\bf C}^1$. If $L$ is a local system (of finite 
dimensional ${\bf C}$ vector spaces) on $U$ we have 

\[ \chi ({\bf P}_{\bf C}^1,j_*L)=2{\rm rk}L-\sum _{i=1}^{p+1}{\rm drop}_{a_i}
(j_*L)\]
with ${\rm drop}_{a_i}(j_*L)={\rm rk}L-{\rm dim}(j_*L)_{a_i}$. The drop 
depends only on the conjugacy class of the local monodromy at $a_i$ and 
decreases under specialization. 

We have a map

\[ {\rm hom} : {\rm conjugacy~classes~in~}GL(n)\times 
{\rm conjugacy~classes~in~}GL(m)\rightarrow {\rm conjugacy~classes~in~}GL(nm)\]

\[ ([A],[B])\mapsto {\rm conjugacy~class~of~}(T\mapsto BTA^{-1})\]
where $T\in {\rm Hom}({\bf C}^n,{\bf C}^m)$. 
If $[A']\subset ~{\rm closure~of~}[A]$, $[B']\subset ~{\rm closure~of~}[B]$, 
then hom$([A'],[B'])\subset$ ${\rm closure~of~hom}([A],[B])$. 

Let $L$ be a local system with local monodromies in $C_j$, $L^*$ a local 
system with local monodromies in $C_j^*$, $L$ being irreducible. By 
assumption $\chi ({\bf P}_{\bf C}^1, j_*\underline{{\rm Hom}}(L,L))=2$. 
Now the local monodromies of $ \underline{{\rm Hom}}(L,L^*)$ have conjugacy 
classes in the closures of the conjugacy classes of the corresponding local 
monodromies of $\underline{{\rm Hom}}(L,L)$, so 

\[ \forall i ~~{\rm drop}_{a_i}\underline{{\rm Hom}}(L,L^*)\leq 
{\rm drop}_{a_i}\underline{{\rm Hom}}(L,L)~,\]
so $\chi ({\bf P}_{\bf C}^1,j_*\underline{{\rm Hom}}(L,L^*))\geq 2$, so 
$H^0(\underline{{\rm Hom}}(L,L^*))\neq 0$ or 
$H^2({\bf P}_{\bf C}^1,j_*\underline{{\rm Hom}}(L,L^*))$ (which is the dual 
of $H^0(\underline{{\rm Hom}}(L,L^*)^{\vee })=
H^0(\underline{{\rm Hom}}(L^*,L))$) is 
non-zero. Since $L$ is irreducible it 
would follow that $L\simeq L^*$ but this contradicts 
the assumption on $C_j^*$.\\ 
 
{\bf Proof of Proposition~\ref{rigidred}:}\\ 

Consider the multiplicative version (the proof in the additive one is 
performed by analogy). Suppose we are 
given the conjugacy classes $C_i\subset GL(n,{\bf C})$, $1\leq i\leq p+1$, 
and we are interested in solutions of 

\begin{equation}\label{star}
M_1\ldots M_{p+1}={\rm id}~,~M_i\in C_i
\end{equation}

We say that a solution $M=(M_1,\ldots ,M_{p+1})$ is {\em rigid} if every 
solution $M'$ in some neighbourhood of $M$ is $GL(n,{\bf C})$-conjugate 
to $M$. Here ``neighbourhood'' can be taken in the classical or in the Zariski 
topology. Excluding the case $n=0$, for a rigid $(p+1)$-tuple the index of 
rigidity is $\geq 1+{\rm dim}({\rm centralizer})\geq 2$ and when the 
centralizer is reduced to scalars ``rigid''$\Leftrightarrow$''index of 
rigidity 
$=2$". We deduce the proposition from the more general one:

\begin{prop}\label{rigidred1}
For a rigid $(p+1)$-tuple, all irreducible subquotients of the representation 
are rigid.
\end{prop}

We say that a finite dimensional representation of $\pi _1$ is 
{\em isotypical} if all its irreducible subquotients are isomorphic. To prove 
Proposition~\ref{rigidred1} we consider first the case of an isotypical 
representation $\rho :\pi _1\rightarrow GL(V)$ of dimension $n=n_1l$ where 
$n_1,l>0$, $\tau :\pi _1\rightarrow GL(W)$ is an $l$-dimensional irreducible 
representation and the semi-simplification of $\rho$ is 
$\tau \oplus \ldots \oplus \tau$ ($n_1$ times). Since the index of rigidity 
increases under specialization we have

\[ (n_1)^2~{\rm index~rig}(\tau )={\rm index~rig}(\tau \oplus \ldots \oplus 
\tau )\geq {\rm index~rig}(\rho )\geq 2~.\]
Since $\tau$ is irreducible, it must have ${\rm index~rig}=2$. 

To do the general case of Proposition~\ref{rigidred1} we suppose that $\rho$ 
is not isotypical. Then one can find an exact sequence of non-zero 
representations

\[ 0\rightarrow V_1\rightarrow V\rightarrow V_2\rightarrow 0\]
where $V_1$ is isotypical and Hom$_{\pi _1}(V_1,V_2)=0$. $V_1$ is the sum of 
all $\tau$-isotypical subrepresentations of $V$ where $\tau$ is an 
irreducible representation of $V$. \\ 

{\bf Claim.} {\em For such an exact sequence if $V$ is rigid, then $V_1$ and 
$V_2$ are rigid. (This will complete the proof by induction.)}\\ 

\begin{lm}
If $\rho ':\pi _1\rightarrow P$ is a representation sufficiently close to 
$\rho$, $GL(V)$-conjugate to $\rho$, then $\rho '$ is $P$-conjugate to $\rho$.
Here $P=\{ g\in GL(V)|gV_1\subset V_1\}$.
\end{lm}

{\em Proof:} 

By general facts on algebraic group actions the conjugacy class of $\rho$ 
is a locally closed subvariety (of the variety of solutions of (\ref{star})) 
isomorphic to $GL(V)/{\rm stabilizer}(\rho )$. 
Hence, $\rho '$ is conjugate to $\rho$ by $g\in GL(V)$ sufficiently close to 
$1$ -- $\rho '=g\rho g^{-1}$, so $g^{-1}(V_1)$ is $\rho$-invariant. 

If $g^{-1}(V_1)\neq V_1$ we get Hom$_{\pi _1}(g^{-1}(V_1),V_2)\neq 0$ since 
the 
projection $V\rightarrow V_2$ defines a non-zero $\pi _1$ map. 
In other words Hom$_{\pi _1}(\rho '|_{V_1},V_2)\neq 0$. If this holds for a 
sequence of $\rho '$'s converging to $\rho$, then since the above is a 
question of solving a homogeneous system of linear equations, i.e. of 
vanishing of certain minors, we get that there is also a non-zero solution 
in the limit case, i.e. Hom$_{\pi _1}(V_1,V_2)\neq 0$ which contradicts the 
assumption.~~~$\Box$\\ 

{\em Proof of the claim:} Suppose that $\rho _1:\pi _1\rightarrow GL(V)$ is 
not rigid. 
Then one has a one-parameter analytic deformation $\rho _{1,\phi }$, 
$\phi \in ({\bf C},0)$, s.t. $\rho _{1,0}=\rho _1$, $\rho _{1,\phi }$ is not 
conjugate to $\rho _{1,0}$ for $\phi \neq 0$. (In general one cannot say that 
the $\rho _{1,\phi }$ for $\phi \neq 0$ are non-equivalent.) Thus given 
analytic deformations of $\rho _1$, $\rho _2$ (within given conjugacy classes 
of local monodromies) it suffices in view of the lemma to extend them to a 
deformation of $\rho$ (within the conjugacy classes).    
  
One can find suitable deformations of $\rho$ of the generators that do not 
necessarily satisfy (\ref{M_j}) and then one tries to correct them by 
conjugation by maps from the parameter space to the unipotent subgroup  

\[ U=\{ \gamma \in GL(V)|(\gamma -1)V\subset V_1,\gamma |_{V_1}=1\} =
\left\{ \left( \begin{array}{cc}1&\ast \\0&1\end{array}\right) \right\} \]
in terms of block matrices. The correction is possible by the following 
analogue of Proposition~\ref{[A_j,X_j]}: 

{\em The map} 

\[ U^{p+1}\rightarrow U~,~(u_1,\ldots ,u_{p+1})\mapsto 
(u_1M_1u_1^{-1})\ldots (u_{p+1}M_{p+1}u_{p+1}^{-1})\] 
{\em has surjective differential.}

We note that $U$ can be identified with 
Lie$(U)=$Hom$_{\bf C}(V_2,V_1)$ and the map above is linear, namely (in 
additive notation for $U$, viewed as representation of $P$)

\[ (u_1,\ldots ,u_{p+1})\mapsto u_1-M_1(u_1)+M_1(u_2-M_2(u_2))+\ldots +
M_1\ldots M_p(u_{p+1}-M_{p+1}(u_{p+1}))~.\]
Thus surjectivity of this map means that the coinvariant space $U_{\pi _1}$ 
is $0$, but its dual is Hom$_{\pi _1}(V_1,V_2)$ whose vanishing is exactly 
the assumption.

\end{document}